\documentclass[12pt]{article}
\usepackage[utf8]{inputenc}
\usepackage{latexsym}     
\usepackage{amsfonts}
\usepackage{amsthm}
\usepackage{amsmath}
\usepackage{amssymb}
\usepackage{graphicx} 
\newtheorem{theorem}{Theorem}[section]

\theoremstyle{remark}

\newtheorem{example}[theorem]{Example}

\def\tm{{\pitchfork}}
\def\mt{{\emptyset}}

\def\im{{\rm Im\ }}

\def\bfe{{\bf e}}

\def\bfk{{\bf k}}
\def\bfi{{\bf i}}
\def\bfj{{\bf j}}

\def\Q{{\mathbb Q}}

\def\T{{\mathbb T}}

\def\d{\partial}

\def\Z{\mathbb{Z}}
\def\R{\mathbb{R}}

\def\T{\mathbb{T}}

\title{Coding for fluids with the transverse intersection algebra\footnote{This work was supported by the US-Israel Binational Science Foundation, Jerusalem}}
\author{Ofir Aharoni, Daniel An\footnote{DA \& AK: SUNY Maritime College}, Alice Kwon, Ruth Lawrence\footnote{OA \& RL: Einstein Institute of Mathematics, Hebrew University of Jerusalem}\\
with an Appendix by Dennis Sullivan\footnote{CUNY Graduate Center NY \& SUNY, Stony Brook NY}}
\date{17.03.25}

\begin{document}

\maketitle
\begin{abstract}

  {\scriptsize
  
   The concept of a fluid algebra was introduced by Sullivan over a decade ago as an algebraic construct which contains everything necessary in order to write down a form of the Euler equation, as an ODE whose solutions have invariant quantities which can be identified as energy and enthalpy. The natural (infinite-dimensional) fluid algebra on co-exact 1-forms on a three-dimensional closed oriented Riemannian manifold leads to an Euler equation which is equivalent to the classical Euler equation which describes non-viscous fluid flow.  In this paper, the recently introduced transverse intersection algebra associated to a cubic lattice of An-Lawrence-Sullivan is used to construct a finite-dimensional fluid algebra on a cubic lattice (with odd periods). The corresponding Euler equation is an ODE which it is proposed is a `good' discretisation of the continuum Euler equation. This paper contains all the explicit details necessary to implement numerically the corresponding Euler equation. Such an implementation has been carried out by our team and results are pending.
      
   }
\end{abstract}
\section{The fluid algebra formulation}
A {\sl fluid algebra} \cite{S10} is a vector space $V$ along with
\begin{enumerate}
    \item[1.] a positive definite inner product $(\ ,\ )$ (the {\sl metric})
    \item[2.] a symmetric non-degenerate bilinear form $\langle\ ,\ \rangle$ (the {\sl linking form})
    \item[3.] an alternating trilinear form $\{\ ,\ ,\ \}$ (the {\sl triple intersection form})
\end{enumerate}
Given a fluid algebra, the {\sl associated Euler equation} is an evolution equation for $X(t)\in{}V$ given implicitly by
\[
(\dot{X},Z)=\{X,DX,Z\}\hbox{ for all test vectors }Z\in{}V\eqno{(1)}
\]
where $D:V\to{}V$ is the operator defined by $\langle{}X,Y\rangle=(DX,Y)$ for all $X,Y\in{}V$.

\medskip\noindent{\bf Invariance of energy $(X,X)$:} \cite{S10}
\[
\frac{d}{dt}(X,X)=2(\dot{X},X)=2\{X,DX,X\}=0
\]

\medskip\noindent{\bf Invariance of helicity $(X,DX)$:} \cite{S10}
\[
\frac{d}{dt}(X,DX)=(\dot{X},DX)+(X,D\dot{X})=2(\dot{X},DX)=2\{X,DX,DX\}=0
\]
where the second equality follows from $(DX,Y)=\langle{}X,Y\rangle=\langle{}Y,X\rangle=(DY,X)=(X,DY)$.

\medskip
\begin{example}
    This is the classical infinite-dimensional example. See \cite{S10}. Let $V$ consist of coexact 1-forms on a 3-dimensional closed oriented Riemannian manifold $M$,
    \[
V=\big\{d^*\omega\big|\omega\in\Omega^2(M)\big\}\subset\Omega^1(M)
    \]
The fluid algebra structures on $V$ are given by
 \[
 (a,b)=\int_Ma\wedge*b,\quad
\langle{}a,b\rangle=\int_Ma\wedge{}db,\quad
\{a,b,c\}=\int_Ma\wedge{}b\wedge{}c
 \]
 Symmetry of $\langle\ ,\ \rangle$ follows from Leibniz $d(a\wedge{}b)=da\wedge{}b-a\wedge{}db$  and Stokes' theorem. Then $D=*d=d^**$. The evolution equation (1) becomes in this case
 \[
\int_M\dot{X}\wedge*Z=\int_MX\wedge*dX\wedge{}Z\hbox{ for all test vectors }Z\in{}V
\]
Rewritten this becomes
\[
\dot{X}=*(X\wedge*dX)+\hbox{ a closed 1-form, }\eta\eqno{(2)}
\]
since those 1-forms $\eta\in\Omega^1(M)$ such that $\int_M\eta\wedge{}d*\omega=0$ for all $\omega\in\Omega^2M$ are those for which $d\eta=0$, by Stokes' theorem. A closed 1-form is (up to homology) an exact 1-form, that is $dp$ for some scalar $p$ (the pressure). For example, for the 3-torus $M=\T^3=\R^3/\Z^3$, a closed 1-form $\eta$ looks like $\eta=dp+adx+bdy+cdz$ where $a,b,c$ are the periods of $\eta$ in the three directions (evaluations of $\int_{\gamma_i}\eta$ where $\gamma_i$, $i=1,2,3$ are loops on the 3-torus in the three directions). Equation (2) is equivalent (up to some considerations on the homology) to the usual formulation of the Euler equation for incompressible fluids,
\[
\d_tu_i+u_j\d_ju_i=\d_ip,\qquad \d_iu_i=0
\]
where $X=u_1dx+u_2dy+u_3dz$.

Recall that an arbitrary 1-form can be written uniquely as a sum of three terms by the Hodge decomposition, an exact 1-form (image of $d$), a coexact 1-form (image of $d^*$) and a harmonic 1-form (intersection of kernels of $d,d^*$). Applying this to $*(X\wedge*dX)$, the sum of the first and third terms is $-\eta$ while the second term is $\dot{X}$. In this way one can write (2) as
\[
\dot{X}=\phi\big(*(X\wedge*dX)\big)\eqno{(3)}
\]
where $\phi:\Omega^1\to{}V$ is the projection defined by the second term of the Hodge decomposition.
\end{example}

\section{The cubical transverse intersection algebra}

In this section we will supply the data for a finite-dimensional fluid algebra which in some way approximates the continuum fluid algebra and thus whose associated Euler equation is hoped to have solutions which exhibit behaviour similar to solutions of the continuum Euler equation. The fluid algebra itself is given in the next section; in this section we define a differential graded algebra on which it is based.

We use a graded subalgebra $FC_*$ of the combinatorial transverse intersection algebra $EC_*$ from \cite{ALS} associated to the three-dimensional cubic lattice, namely the part in dimension at most two generated by squares of side length $2h$. All these algebras are considered over the field $\Q$.
To make this paper self-contained, we here give an explicit description of $FC_*$, without reference to $EC_*$.

Consider a periodic three-dimensional lattice $\Lambda$ with period $N$ in each direction; denote the lattice spacing by $h$. The dimension zero part $FC_0$ will have basis $\{\mt_a|a\in\Lambda\}$ indexed by vertices in the lattice, $\dim{FC_0}=N^3$. The dimension two part $FC_2$ will have basis $\bigcup_{a\in\Lambda}\{yz_a,xz_a,xy_a\}$ indexed by squares of side $2h$ in the lattice, parallel to one of the coordinate planes, $\dim{FC_2}=3N^3$ with $a$ specifying the centre of a square. The dimension one part $FC_1$ will have basis indexed by sticks $x_{ab}, y_{ab}, z_{ab}$ of length $h$ in the lattice and infinitesimal sticks $x_{aa},y_{aa},z_{aa}$ in the lattice (with direction parallel to one of the coordinate axes), $\dim{FC_1}=6N^3$; here $a,b\in\Lambda$ label the endpoints of the sticks where their difference $b-a\in\{h\bfe_1,h\bfe_2,h\bfe_3\}$. Grading is by codimension (three minus the dimension).

The transverse multiplication on $FC_*$ is induced by that on $EC_*$ in \cite{ALS} where each $2h$-square in $FC_2$ is considered as a sum of the four $h$-squares generators of $EC_2$ into which it geometrically crumbles. Explicitly the non-zero products are $FC_1\times{}FC_2\to{}FC_0$ 
and $FC_2\times{}FC_2\to{}FC_1$. 

For $\tm:FC_1\times{}FC_2\to{}FC_0$, the non-zero products of basis vectors are of orthogonal intersecting sticks ($h$- or infinitesimal) and squares, where the product is a signed multiple of the intersecting point; those multiples depend on the type and relative location of the intersection and are given here
\[\includegraphics[width=.7\textwidth]{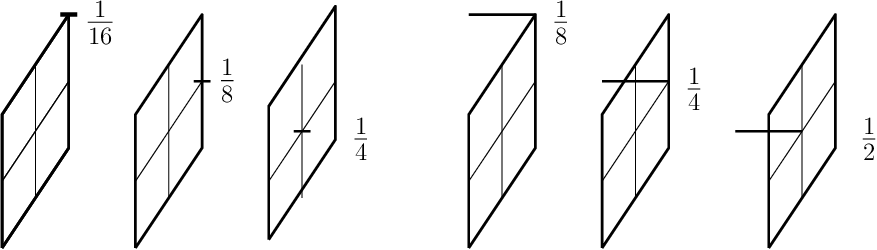}\]
The sign is given by the relative orientation of the square and interval.
These all follow from the products in $EC_*$,
\[\includegraphics[width=.1\textwidth]{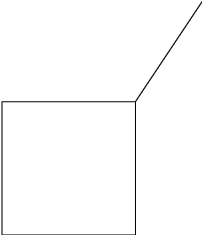}\hskip1ex\raise5ex\hbox{product$=\frac18\cdot$}\qquad
\includegraphics[width=.04\textwidth]{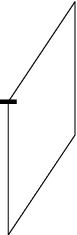}\hskip1em\raise5ex\hbox{product$=\frac{1}{16}\cdot$}\]

For $\tm:FC_2\times{}FC_2\to{}FC_1$, the non-zero products of basis elements involve $2h$-squares in orthogonal directions which intersect geometrically, in a point, an $h$-stick or a $2h$-stick. There are several possible configurations. For a  point intersection, the intersections are all multiples of an infinitesimal element in the common direction
\[\includegraphics[width=.7\textwidth]{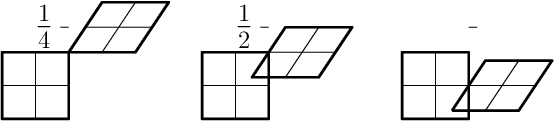}\]
For an $h$-stick intersection, the intersections are multiples of the $h$-stick given by the geometric intersection
\[\includegraphics[width=.7\textwidth]{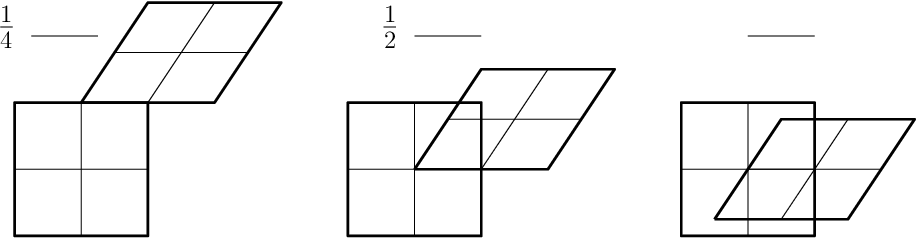}\]
For a $2h$-stick intersection, all the intersections are muliples of a $2h$-stick (that is, a sum of two $h$-sticks) minus infinitesimal sticks at the endpoints of the $2h$-stick geometric intersection,
\[\includegraphics[width=.7\textwidth]{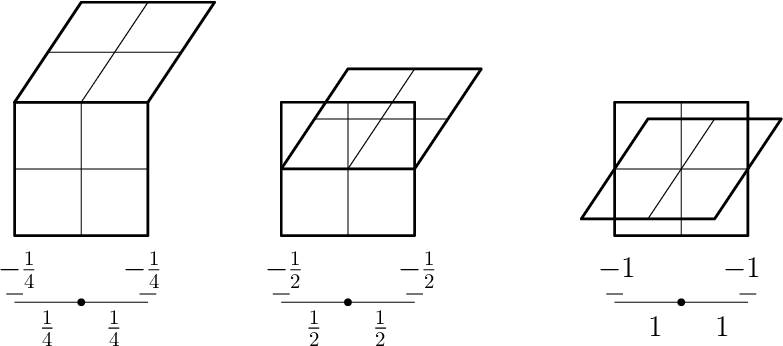}\]

The boundary operator $\d$ is defined in the usual way, $\d\mt_a=0$, $\d{}u_{ab}=\mt_b-\mt_a$ for any symbol $u\in\{x,y,z\}$, while the boundary of a $2h$-square is the signed sum of the four $2h$-edges, that is of eight $h$-sticks. 

The properties of $\tm$ and $\d$ on $FC_*$ from \cite{ALS} are that they define a differential graded commutative associative algebra. Here the product rule
\[
\d(a\tm{}b)=(\d{}a)\tm{}b+(-1)^{c(a)}a\tm\d{}b
\]
holds for any $a,b\in{}FC_*$, the only non-trivial case being $a,b\in{}FC_2$. Furthermore, the natural crumbling map $FC_*\to{}FC'_*$ from the transverse intersection algebra of a coarse cubical lattice to the transverse intersection algebra of a finer cubical lattice commutes with both $\d$ and $\tm$.

\section{Fluid algebra induced by $FC_*$ for odd $N$}
Using the cubic transverse intersection algebra of the previous section, we here define a finite-dimensional fluid algebra on a vector space $V\subset{}FC_2$ in the case of $N$ odd.
We start by clarifying some of the geometry and linear algebra in $FC_*$.

According to the type of basis element, infinitesimal stick or $h$-stick, $FC_1$ can be decomposed as 
\[
FC_1=I_1\oplus{}C_1
\]
where $I_1$ is spanned by infinitesimal sticks and $C_1$ is spanned by $h$-sticks (both have dimension $3N^3$). Considering a $2h$-stick as a sum of two $h$-sticks, thus $2h$-sticks are elements of $C_1$; since $N$ is odd, $2h$-sticks also form a basis for $C_1$.

Next observe that there is a natural map $*:FC_2\to{}C_1$ which maps a $2h$-square to the orthogonal $2h$-stick sharing the centre point with that of the square (with sign matching the orientation). And similarly the inverse map $*:C_1\to{}FC_2$ maps a $2h$-stick to the orthogonal $2h$-square sharing the centre point with that of the stick.

Let $\#:FC_0\to\Q$ be the augmentation map which counts points (with weighting) that is $\mt_a\mapsto1$.

Analogously to the continuum example, the vector space of the fluid algebra is defined by 
\[
V=\big\{*\d{}x\big|x\in{}FC_2\big\}
\]
The structure maps in the fluid algebra are given by
\[
(a,b)\equiv\#(a\tm*b),\qquad\langle{}a,b\rangle\equiv\#(a\tm\d{}b),\qquad
\{a,b,c\}\equiv\#(a\tm{}b\tm{}c)
\]

\noindent{\bf Verification of fluid algebra properties:}

\noindent{\sl Metric:} The metric is positive definite on all of $FC_2$; indeed the inner product defined by $(\ ,\ )$ is non-zero on pairs of basis elements of $EC_2$ only when they are parallel and have centres at most $h$ apart (by $\|\ \|_\infty$) while their inner product is the power of $\frac12$ given by the $d_1$ norm of their centre separation. The matrix of inner products is thus three diagonal copies (for the three possible square orientations) of the third tensor-power of the equivalent one-dimensional matrix, a symmetric circular matrix with $1$ on the diagonal and $\frac12$ on the off-diagonals (including the opposite corner elements). Its eigenvalues are $1+\frac12\omega^j+\frac12\omega^{-j}$ for $\omega=\exp{\frac{2\pi{}i}{N}}$, $j=0,1,\ldots,N-1$; these simplify to $1+\cos{\frac{2\pi{}j}{N}}$ which are all positive for $N$ odd. This verifies positive definiteness on $FC_2$ and hence also on the subspace $V$.

\medskip\noindent{\sl Linking form:} Symmetry of $\langle\ ,\ \rangle$ follows from $\d(a\tm{}b)=b\tm\d{}a-a\tm\d{}b$ along with the fact that $\#(\d{}x)=0$ for any $x\in{}FC_1$. Note that $V$ as defined, is the star of the space of 1-chains spanned by boundaries of 2-chains. There are $3N^3$ boundaries of basis 2-chains, amongst which there are $3N$ independent relations so that $\dim{V}=3(N^3-N)$. 

To prove non-degeneracy of the linking form on $V$, suppose that $a\in{}V$ is such that $\langle{}a,b\rangle=0$ for all $b\in{}V$. We wish to prove that $a=0$. Choose $b=*\d{}a$. By definition $b\in{}V$ and so
\[
0=\langle{}a,*\d{}a\rangle=\#(\d{}a\tm*\d{}a)=(*\d{a},*\d{}a)
\]
from which it follows that $*\d{}a=0$ since the metric is positive definite (and thus the inner product non-degenerate). Since $a\in{}V$, it can be written in the form $a=*\d{}x$ for some $x\in{}FC_2$. Thus we now have $\d*\d{}x=0$. Hence
\[
0=\#(x\tm\d*\d{}x)=\#(\d{}x\tm*\d{}x)=(*\d{}x,*\d{}x)=(a,a)
\]
for which it follows by non-degeneracy of the metric that $a=0$ as required.

\medskip\noindent{\sl Triple bracket:} The alternating property of $\{\ ,\ ,\ \}$ follows from associativity of $\tm$ along with graded commutativity. Since the grading is by codimension, this is odd for $EC_2$ and hence such elements anti-commute.

\medskip\noindent{\bf Associated Euler equation:} By definition $D=*\d$; this is meaningful as a map $FC_2\to{}FC_2$ or even $V\to{}V$ since $\d:FC_2\to{}C_1$. The evolution equation (1) for the 2-chain $X(t)$ becomes
\[
\#(*\dot{X}\tm{}Z)=\#(X\tm*\d{}X\tm{}Z)\hbox{ for all }Z\in{}V
\]
That is, $*\dot{X}$ is an element of $\im{\d}\subset{}C_1$ whose intersection number with every 2-chain in $V$ is identical with that of $X\tm*\d{}X$. Since the metric is non-degenerate, this specifies $\dot{X}\in{}V$ uniquely. More precisely, there is a unique linear map $\phi:FC_1\to{}V$ for which
\[
\#(*\phi(a)\tm{}v)=\#(a\tm{}v)\hbox{ for all }v\in{}V,\ a\in{}FC_1
\]
Observing the definition of the intersection $\tm:FC_1\times{}FC_2\to{}FC_0$ in the previous section, it may be seen that an $h$-stick has exactly the same intersection numbers with elements of $FC_2$ as a sum of two infinitesimal sticks at its endpoints, each with weighting 2. Use this to define a map $i:FC_1\to{}I_1$ for which $\#(a\tm{}v)=\#(i(a)\tm{}v)$ for all $a\in{}FC_1$ and $v\in{}FC_2$. Then $\phi$ factors through $i$; that is, there is a map $\pi:I_1\to{}V$ such that $\phi=\pi\circ{}i$. Using this notation, Euler's equation in the fluid algebra can be written as
\[
\dot{X}=\pi\big(i(X\tm*\d{}X)\big)\eqno{(4)}
\]
The coefficients of $i(X\tm*\d{}X)$ are given by Dan's formula of 428 terms or equivalently my formula as a sum of eighteen products of sums of three terms with sums of ten terms. See the computation below.

It remains to compute $\pi$ (``solve the Poisson problem")  which is the unique map $\pi:I_1\to{}V$ defined by 
\[
\#(*\pi(a)\tm{}v)=\#(a\tm{}v)\hbox{ for all }v\in{}V,\ a\in{}I_1
\]
This we will accomplish in the next section.

\medskip\noindent{\bf Computation of $i(X\tm*\d{}X)$:} Here $X\in{}FC_2$ is a linear combination of $2h$-squares,
\[
X=\sum\limits_{a\in\Lambda}u^{yz}_ayz_a+u^{zx}_azx_a+u^{xy}_axy_a
\]
where $u^{yz}_a,u^{zx}_a,u^{xy}_a$ are the coefficients of the various square types. 
The boundary of a $2h$-square is given as a sum of four $2h$-sticks, for example
\[
\d(yz_a)=y_{a-h\bfk}+z_{a+h\bfj}-y_{a+h\bfk}-z_{a-h\bfj}
\]
and so its star is a combination of four $2h$-squares,
\begin{align*}
    *\d(yz_a)&=zx_{a-h\bfk}+xy_{a+h\bfj}-zx_{a+h\bfk}-xy_{a-h\bfj}\\
    *\d(zx_a)&=xy_{a-h\bfi}+yz_{a+h\bfk}-xy_{a+h\bfi}-yz_{a-h\bfk}\\
    *\d(xy_a)&=yz_{a-h\bfj}+zx_{a+h\bfi}-yz_{a+h\bfj}-zx_{a-h\bfi}
\end{align*}
Putting this together
\begin{align*}
*\d{}X=&\sum\limits_{a\in\Lambda}(u^{zx}_{a-h\bfk}-u^{zx}_{a+h\bfk}+u^{xy}_{a+h\bfj}-u^{xy}_{a-h\bfj})yz_a\\
&+(u^{xy}_{a-h\bfi}-u^{xy}_{a+h\bfi}+u^{yz}_{a+h\bfk}-u^{yz}_{a-h\bfk})zx_a
+(u^{yz}_{a-h\bfj}-u^{yz}_{a+h\bfj}+u^{zx}_{a+h\bfi}-u^{zx}_{a-h\bfi})xy_a
\end{align*}
which we will write as 
\[
*\d{}X=\sum\limits_{a\in\Lambda}v^{yz}_ayz_a+v^{zx}_azx_a+v^{xy}_axy_a
\]
where $v^{yz}_a=u^{zx}_{a-h\bfk}-u^{zx}_{a+h\bfk}+u^{xy}_{a+h\bfj}-u^{xy}_{a-h\bfj}$,
$v^{zx}_a=u^{xy}_{a-h\bfi}-u^{xy}_{a+h\bfi}+u^{yz}_{a+h\bfk}-u^{yz}_{a-h\bfk}$ and
$v^{xy}_a=u^{yz}_{a-h\bfj}-u^{yz}_{a+h\bfj}+u^{zx}_{a+h\bfi}-u^{zx}_{a-h\bfi}$. Now take the transverse intersection product with $X$ and get
\begin{align*}
X\tm*\d{}X=&\sum\limits_{a,b\in\Lambda}(u_a^{zx}v_b^{xy}-u_b^{xy}v_a^{zx})zx_a\tm{}xy_b\\
&+(u_a^{xy}v_b^{yz}-u_b^{yz}v_a^{xy})xy_a\tm{}yz_b+(u_a^{yz}v_b^{zx}-u_b^{zx}v_a^{yz})yz_a\tm{}zx_b
\end{align*}
The transverse intersection product $zx_a\tm{}xy_b$ yields a one-dimensional object in the $x$-direction, so long as there is a non-empty geometric intersection of the associated $2h$-squares, that is so long as $|a_1-b_1|\leq2h$, $|a_2-b_2|,|a_3-b_3|\leq{}h$. In this case, the intersection will be $(\frac12)^{\frac1h(|a_2-b_2|+|a_3-b_3|)}$ times
\begin{align*}
    (x)_{\frac{a_1+b_1}{2},a_2,b_3}&\hbox{ when $|a_1-b_1|=2h$}\\
    x_{[a_1,b_1],a_2,b_3}&\hbox{ when $|a_1-b_1|=h$}\\
    x_{a_1,a_2,b_3}-(x)_{a_1-h,a_2,b_3}-(x)_{a_1+h,a_2,b_3}&\hbox{ when $a_1=b_1$}
\end{align*}
where $(x)_a$ denotes the infinitesimal edge in direction $x$ located at $a$, $x_a$ denotes the $2h$-stick in direction $x$ whose centre is at $a$ while $x_{[a_1,b_1],a_2,b_3}$ denotes the $h$-stick in direction $x$ between the points $(a_1,a_2,b_3)$ and $(b_1,a_2,b_3)$. Recall that a $2h$-stick in $FC_1$ is equal to the sum of the two $h$-sticks into which it decomposes while the map $i:FC_1\to{}I_1$ is defined to be the identity on infinitesimal sticks and to replace $h$-sticks by an equivalent combination of infinitesimal sticks (equivalent in the sense that all transverse intersections with $2h$-squares remain unaltered). Thus $i$ acts on $x$-directed sticks by
\[
(x)_a\mapsto(x)_a,\qquad{}x_{ab}\mapsto2(x)_a+2(x)_b,\qquad{}x_a\mapsto2(x)_{a-h\bfi}+4(x)_a+2(x)_{a+h\bfi}
\]
and hence $i(zx_a\tm{}xy_b)$ when non-zero is $(\frac12)^{\frac1h(|a_2-b_2|+|a_3-b_3|)}$ times
\begin{align*}
    (x)_{\frac{a_1+b_1}{2},a_2,b_3}&\hbox{ when $|a_1-b_1|=2h$}\\
    2(x)_{a_1,a_2,b_3}+2(x)_{b_1,a_2,b_3}&\hbox{ when $|a_1-b_1|=h$}\\
    (x)_{a_1-h,a_2,b_3}+4(x)_{a_1,a_2,b_3}+(x)_{a_1+h,a_2,b_3}&\hbox{ when $a_1=b_1$}
\end{align*}
The coefficient of $(x)_c$ in $i(X\tm\d{}X)$ is thus 
\[
\sum(\frac12)^{\frac1h(|a_2-b_2|+|a_3-b_3|)}(u_a^{zx}v_b^{xy}-u_b^{xy}v_a^{zx})(\hbox{1, 2 or 4})
\]
over those $a,b\in\Lambda$ with $a_2=c_2$, $b_3=c_3$, $|a_2-b_2|\leq{}h$, $|a_3-b_3|\leq{}h$, while
\[
c_1=\left\{
\begin{array}{ll}
(a_1+b_1)/2&\hbox{ when $|a_1-b_1|=2h$}\\
a_1 \hbox{ or }b_1&\hbox{ when $|a_1-b_1|=h$}\\
a_1-h\hbox{ or }a_1\hbox{ or }a_1+h&\hbox{ when $a_1=b_1$}
\end{array}\right.
\]
and the term is counted with factor 1, 2 or 4 according to the coefficients in the above formula for $i(zx_a\tm{}xy_b)$. This give three choices for $b_2$, namely $c_2,c_2\pm{}h$ (weights $1,\frac12$) and similarly three choices for $a_3$, namely $c_3,c_3\pm{}h$ (weights $1,\frac12$) while there are nine choices for $(a_1,b_1)$,
\[
(a_1,b_1)=(c_1-h,c_1+h),(c_1+h,c_1-h),(c_1,c_1\pm{}h),(c_1\pm{}h,c_1),(c_1,c_1),(c_1\pm{}h,c_1\pm{}h)
\]
with weights $1,1,2,2,4,1$ respectively. Finally we get the following formula for the coefficient of $(x)_c$ in $i(X\tm\d{}X)$,
\begin{align*}
&\sum\limits_{(a_1,b_1)}(\hbox{1, 2 or 4})(U^{zx,z}_{a_1,c_2,c_3}V_{b_1,c_2,c_3}^{xy,y}-U_{b_1,c_2,c_3}^{xy,y}V_{a_1,c_2,c_3}^{zx,z})\\
&=U^{zx,z}_{c-h\bfi}V_{c+h\bfi}^{xy,y}-U_{c+h\bfi}^{xy,y}V_{c-h\bfi}^{zx,z}
+U^{zx,z}_{c+h\bfi}V_{c-h\bfi}^{xy,y}-U_{c-h\bfi}^{xy,y}V_{c+h\bfi}^{zx,z}\\
&+2U^{zx,z}_{c}V_{c+h\bfi}^{xy,y}-2U_{c+h\bfi}^{xy,y}V_{c}^{zx,z}
+2U^{zx,z}_{c}V_{c-h\bfi}^{xy,y}-2U_{c-h\bfi}^{xy,y}V_{c}^{zx,z}\\
&+2U^{zx,z}_{c+h\bfi}V_{c}^{xy,y}-2U_{c}^{xy,y}V_{c+h\bfi}^{zx,z}
+2U^{zx,z}_{c-h\bfi}V_{c}^{xy,y}-2U_{c}^{xy,y}V_{c-h\bfi}^{zx,z}\\
&+4U^{zx,z}_{c}V_{c}^{xy,y}-4U_{c}^{xy,y}V_{c}^{zx,z}
+U^{zx,z}_{c+h\bfi}V_{c+h\bfi}^{xy,y}-U_{c+h\bfi}^{xy,y}V_{c+h\bfi}^{zx,z}
+U^{zx,z}_{c-h\bfi}V_{c-h\bfi}^{xy,y}-U_{c-h\bfi}^{xy,y}V_{c-h\bfi}^{zx,z}
\end{align*}
where $U$ and $V$ are smeared versions of $u$,$v$ in the directions indicated by the indices, that is
\begin{align*}
    U^{xy,x}_a=&u^{xy}_a+\frac12u^{xy}_{a-h\bfi}+\frac12u^{xy}_{a+h\bfi}\\
    U^{xy,y}_a=&u^{xy}_a+\frac12u^{xy}_{a-h\bfj}+\frac12u^{xy}_{a+h\bfj}\\
    U^{zx,z}_a=&u^{zx}_a+\frac12u^{zx}_{a-h\bfk}+\frac12u^{zx}_{a+h\bfk}\\
    U^{zx,x}_a=&u^{zx}_a+\frac12u^{zx}_{a-h\bfi}+\frac12u^{zx}_{a+h\bfi}\\
    U^{yz,y}_a=&u^{yz}_a+\frac12u^{yz}_{a-h\bfj}+\frac12u^{yz}_{a+h\bfj}\\
    U^{yz,z}_a=&u^{yz}_a+\frac12u^{yz}_{a-h\bfk}+\frac12u^{yz}_{a+h\bfk}
\end{align*}
while for example $V^{xy,x}_a=v^{xy}_a+\frac12v^{xy}_{a-h\bfi}+\frac12v^{xy}_{a+h\bfi}$
so that
\begin{align*}
V^{xy,x}_a=&u^{yz}_{a-h\bfj}-u^{yz}_{a+h\bfj}+u^{zx}_{a+h\bfi}-u^{zx}_{a-h\bfi}\\
    &+\frac12(u^{yz}_{a-h\bfi-h\bfj}-u^{yz}_{a-h\bfi+h\bfj}+u^{yz}_{a+h\bfi-h\bfj}-u^{yz}_{a+h\bfi+h\bfj}+u^{zx}_{a+2h\bfi}-u^{zx}_{a-2h\bfi})\\
V^{xy,y}_a=&u^{yz}_{a-h\bfj}-u^{yz}_{a+h\bfj}+u^{zx}_{a+h\bfi}-u^{zx}_{a-h\bfi}\\
    &+\frac12(u^{zx}_{a+h\bfi+h\bfj}-u^{zx}_{a-h\bfi+h\bfj}+u^{zx}_{a+h\bfi-h\bfj}-u^{zx}_{a-h\bfi-h\bfj}+u^{yz}_{a-2h\bfj}-u^{yz}_{a+2h\bfj})\\
V^{zx,z}_a=&u^{xy}_{a-h\bfi}-u^{xy}_{a+h\bfi}+u^{yz}_{a+h\bfk}-u^{yz}_{a-h\bfk}\\
    &+\frac12(u^{xy}_{a-h\bfi-h\bfk}-u^{xy}_{a+h\bfi-h\bfk}+u^{xy}_{a-h\bfi+h\bfk}-u^{xy}_{a+h\bfi+h\bfk}+u^{yz}_{a+2h\bfk}-u^{yz}_{a-2h\bfk})\\
V^{zx,x}_a=&u^{xy}_{a-h\bfi}-u^{xy}_{a+h\bfi}+u^{yz}_{a+h\bfk}-u^{yz}_{a-h\bfk}\\
    &+\frac12(u^{yz}_{a+h\bfi+h\bfk}-u^{yz}_{a+h\bfi-h\bfk}+u^{yz}_{a-h\bfi+h\bfk}-u^{yz}_{a-h\bfi-h\bfk}+u^{xy}_{a-2h\bfi}-u^{xy}_{a+2h\bfi})\\
V^{yz,y}_a=&u^{zx}_{a-h\bfj}-u^{zx}_{a+h\bfj}+u^{xy}_{a+h\bfi}-u^{xy}_{a-h\bfi}\\
    &+\frac12(u^{zx}_{a-h\bfj-h\bfk}-u^{zx}_{a-h\bfj+h\bfk}+u^{zx}_{a+h\bfj-h\bfk}-u^{zx}_{a+h\bfj+h\bfk}+u^{xy}_{a+2h\bfj}-u^{xy}_{a-2h\bfj})\\
V^{yz,z}_a=&u^{zx}_{a-h\bfk}-u^{zx}_{a+h\bfk}+u^{xy}_{a+h\bfj}-u^{xy}_{a-h\bfj}\\
    &+\frac12(u^{xy}_{a+h\bfj+h\bfk}-u^{xy}_{a-h\bfj+h\bfk}+u^{xy}_{a+h\bfj-h\bfk}-u^{xy}_{a-h\bfj-h\bfk}+u^{zx}_{a-2h\bfk}-u^{zx}_{a+2h\bfk})
\end{align*}
Notice that a cancellation occurs when evaluating $V$ leaving it as a sum of ten terms rather than the expected twelve; for example $u^{zx}_a$ appears in $V_a^{xy,x}$ with both a negative and positive sign.
The final formula for the coefficients in $i(X\tm\d{}X)$ is a sum of 18 terms, each of which is a product of a sum of three evaluations of $u$ with a sum of ten evaluations of $u$; in total 540 terms of which a number cancel leaving Dan's formula of a sum of 428 terms.

\section{Solving the Poisson problem}
The only part of the evolution equation (4) which we have not yet given explicitly is the map $\pi:I_1\to{}V$, uniquely defined by the property
\[
\#(*\pi(a)\tm{}v)=\#(a\tm{}v)\hbox{ for all }v\in{}V,\ a\in{}I_1
\]
In particular we will find in this section a formula for the image under $\pi$ of a basis element of $I_1$, the infinitesimal stick $(x)_0$.

\medskip First we define a map $r:I_1\to{}C_1$ which is such that $r(a)\tm{}v=a\tm{}v$ for all $v\in{}FC_2$, $a\in{}I_1$. Recall that up to transverse intersections with elements of $FC_2$, an $h$-stick is equivalent to a linear combination of infinitesimal sticks,
\[
x_{a,a+h\bfi}\sim2(x)_a+2(x)_{a+h\bfi}
\]
Inverting this relation, we find that 
\[
(x)_0\sim\sum_{k=0}^{N-1}\frac14(-1)^kx_{kh\bfi,(k+1)h\bfi}\eqno{(5)}
\]
once again the sign $\sim$ denoting an equivalence with regard to transverse intersection with squares in $FC_2$. Up to symmetry and translation, this defines the map $r:I_1\to{}C_1$. Recall that $C_1$ has a basis consisting of $2h$-sticks; with respect to this basis
\[
r((x)_0)=\sum\limits_{k=0}^{N-1}(-1)^k\frac{N-2k}{8}x_{r\bfi}\eqno{(6)}
\]

Consider now the complex $C_*$ which is generated by cubic cells of edge length $2h$ and dimensions 0,1,2,3  on the lattice $\Lambda$. Because $N$ is odd, these chain spaces equivalently have bases given by cubic cells of side length $h$ in the matching dimensions; nonetheless we prefer to use the $2h$-basis because it affords a natural description for the star operator. This complex becomes a chain complex under the geometric boundary $\d:C_j\to{}C_{j-1}$. There is a natural star operator $*:C_j\to{}C_{3-j}$ which takes a $2h$-cell to the (signed) complementary cell sharing the same centre. The operator $\d^*=*\d*:C_j\to{}C_{j+1}$ is adjoint (up to sign) to $\d$ with respect to the symmetric non-degenerate pairing $(x,y)=\#(x\tm*y)$ as 
\begin{align*}
(\d{}a,b)=\#(\d{}a\tm*b)&=\#\big(\d(a\tm*b)-(-1)^{c(a)}a\tm\d*b\big)\\
&=(-1)^{|a|}\#(a\tm*\d^*b)=(-1)^{|a|}(a,\d^*b)
\end{align*}
By the Hodge decomposition, any $q\in{}C_1$ can be decomposed uniquely as a sum
\[
q=e+f+c\eqno{(7)}
\]
with $e\in\im(\d:C_2\to{}C_1)$, $f\in\im(\d^*:C_0\to{}C_1)$, $c\in\ker\d\cap\ker\d^*$. The terms $e,f,c$ are all linear functions $e(q),f(q),c(q)$ of $q$. Note that the harmonic elements of $C_1$, that is those in $\ker(\d:C_1\to{}C_0)\cap\ker(\d^*:C_1\to{}C_2)$ are  constant combinations
\[c=\sum\limits_{a\in\Lambda}c_1x_a+c_2y_a+c_3z_a\]
for some constants $c_1,c_2,c_3$ independent of lattice point $a\in\Lambda$. The harmonic element $c$ determined from $q\in{}C_1$ by (7) is given by choosing these constants $c_1,c_2,c_3$ to be the averages over all lattice points of the coefficients of $2h$-sticks in the three directions in $q$. 

The map $\pi:I_1\to{}V$ is given by $\pi(g)=*e(r(g))\in\im(*\d)=V$. To verify the defining property, note from (7) that $\d^*q=\d^*e$ so that $\d*q=\d*e$. Take any $v\in{}V$, say $v=*\d{}w$, $w\in{}FC_2$
\begin{align*}
\#(*\pi(g)\tm{}v)&=\#(e(r(g))\tm*\d{}w)=\#(\d*e(r(g))\tm{}w)\\
&=\#\big((\d*{}r(g))\tm{}w\big)=\#(r(g)\tm*\d{}w)=\#(r(g)\tm{}v)=\#(g\tm{}v)
\end{align*}

So in order to find the map $e$ in (7), it suffices to write it as $e(q)=q-f(q)-c(q)$ and find the map $f(q)$. That is,
\[
*\pi(g)=r(g)-c(r(g))-f(r(g))\eqno{(8)}
\]
We already gave the formula for the constant part $c(q)$ above in terms of averages in the three directions. From (7), $\d{}q=\d{}f$ while $f\in\im(\d^*)$ so that $f=\d^*p$ for some $p\in{}C_0$. To find $f$ it is sufficient to solve Poisson's equation
\[\d\d^*p=\d{q}\eqno{(9)}\]
which has a unique solution for $p\in{}C_0$ (a function on the lattice) up to the addition of a constant. 

We now carry out this process on an infinitesimal $x$-stick to find $\pi((x)_0)\in{}V$. First of all $q=r((x)_0)\in{}C_1$ is given by (5), (6). The constant part $c(q)$ is given by taking the average coefficients in the three directions; the sum of coefficients in the $x$-direction on $2h$-sticks in (6) is $\frac18$ and so 
\[
c(q)=\frac1{8N^3}\sum\limits_{a\in\Lambda}x_a\eqno{(10)}
\]
Next $f(q)=\d^*p$ where $p\in{}C_0$ is a solution of Poisson's equation (9). Taking the boundary of (5) we get
\[
\d{}q=\sum\limits_{k=1}^{N-1}\frac12(-1)^{k-1}\mt_{kh\bfi}
\]
Writing (9) as an equation for the scalar function $p_a$, $a\in\Lambda$, we have
\[p_{a+2h\bfi}+p_{a-2h\bfi}+p_{a+2h\bfj}+p_{a-2h\bfj}+p_{a+2h\bfk}
+p_{a-2h\bfk}-6p_{a}\\
=\frac12(-1)^{a_1-1}\delta_{a_1\not=0}\delta_{a_2=0}\delta_{a_3=0}\eqno{(11)}
\]
where the left hand side is the Laplacian $\Delta{}p$.

\bigskip\noindent{\bf Solving for $p$:} Our technique to solve (11) is to find Laplacians of a selection of functions a suitable linear combination of which will provide the solution. Rather than considering $p$ as a function on the periodic lattice (torus) with spacing $h$, we will label the lattice points by $(x,y,z)$ coordinates each running amongst $0,1,\ldots,N-1$ (equivalently take $h=1$) and consider the values at these $N^3$ points as determining all values by periodic extension. Consider a function
\[
p_{i,j,k}=\lambda^i\mu^j\nu^k
\]
If this function was on the infinite lattice $\Z^3$, it would be harmonic under the condition
\[\lambda^2+\lambda^{-2}+\mu^2+\mu^{-2}+\nu^2+\nu^{-2}=6\eqno{(12)}\]
If $\lambda^N=\mu^N=1$ then this function is periodic in the $x$ and $y$ directions. Truncating this function to $B=\{0,1,\ldots,N-1\}^3$ and then extending periodically we find that $\Delta{}p$ vanishes everywhere except close to the `discontinuity' over $z=0$; to be precise, it vanishes everywhere on $B$ at except at points with $z$-coordinate $0,1,N-2,N-1$ where
\[
(\Delta{}p)_a=\left\{
\begin{array}{ll}
\lambda^{a_1}\mu^{a_2}\nu^{a_3}(\nu^{N-2}-\nu^{-2}),&\hbox{ if $a_3=0,1$}\\
\lambda^{a_1}\mu^{a_2}\nu^{a_3}(\nu^{2-N}-\nu^{2}),&\hbox{ if $a_3=N-2,N-1$}
\end{array}
\right.
\]
Setting $\omega=\exp{\frac{2\pi{}i}{N}}$, an $N$-th root of unity, we have for each $\lambda=\omega^l$, $\mu=\omega^m$ with $l,m\in\{0,1,\ldots,N-1\}$ four real solutions $\nu$ to (12),
\[\nu^2+\nu^{-2}=6-2\cos\frac{4\pi{}l}{N}-2\cos\frac{4\pi{}m}{N}\eqno{(13)}\]
which we will write as $\nu=\pm\nu_{l,m},\pm\nu_{l,m}^{-1}$ where $\nu_{l,m}$ is the largest solution. Explicitly the roots are
\[\nu=\pm\sqrt{1-\frac12\cos\frac{4\pi{}l}{N}-\frac12\cos\frac{4\pi{}m}{N}}
\pm\sqrt{2-\frac12\cos\frac{4\pi{}l}{N}-\frac12\cos\frac{4\pi{}m}{N}}\]
while the largest is
\[\nu_{l,m}=\sqrt{1-\frac12\cos\frac{4\pi{}l}{N}-\frac12\cos\frac{4\pi{}m}{N}}
+\sqrt{2-\frac12\cos\frac{4\pi{}l}{N}-\frac12\cos\frac{4\pi{}m}{N}}\eqno{(14)}\]
These solutions will be distinct except in the case $l=m=0$ when $\nu_{0,0}=1$.

Looking at the Laplacians, $\Delta{}p$, of these four solutions, for fixed $\lambda=\omega^l$, $\mu=\omega^m$ but the four different values of $\nu$, we can take a suitable linear combination so that the Laplacian vanishes everywhere but on the plane $a_3=0$. Indeed the Laplacian of each of the four solutions vanishes everywhere except on the four planes $a_3=0,1,N-2,N-1$ and there it is a multiple of $\lambda^{a_1}\mu^{a_2}$, the multiple being 
\begin{align*}
    \nu:\quad&\nu^{N-2}-\nu^{-2},\nu^{N-1}-\nu^{-1},1-\nu^N,\nu-\nu^{N+1}\\
   -\nu:\quad&-\nu^{N-2}-\nu^{-2},\nu^{N-1}+\nu^{-1},1+\nu^{N},-\nu-\nu^{N+1}\\
   \nu^{-1}:\quad&\nu^{2-N}-\nu^{2},\nu^{1-N}-\nu,1-\nu^{-N},\nu-\nu^{-N-1}\\
    -\nu^{-1}:\quad&-\nu^{2-N}-\nu^{2},\nu^{1-N}+\nu,1+\nu^{-N},-\nu-\nu^{-N-1}
\end{align*}
Combining the first two lines, for $p_a=\frac12\lambda^{a_1}\mu^{a_2}
\left(\frac{\nu^{a_3}}{\nu^{N-2}-\nu^{-2}}-\frac{(-\nu)^{a_3}}{\nu^{N-2}+\nu^{-2}}\right)$ we see that $(\Delta{}p)_a$ vanishes except when $a_3=0,N-2$ and there it is $1,-\nu^2$ respectively times $\lambda^{a_1}\mu^{a_2}$. Combining with the third and fourth lines gives that
\[
p_a=\frac12\frac{\lambda^{a_1}\mu^{a_2}}{\nu^2-\nu^{-2}}
\left(\frac{\nu^{a_3}}{1-\nu^{N}}+\frac{(-\nu)^{a_3}}{1+\nu^N}
-\frac{\nu^{-a_3}}{1-\nu^{-N}}-\frac{(-\nu)^{-a_3}}{1+\nu^{-N}}
\right)
\]
has $(\Delta{}p)_a=\lambda^{a_1}\mu^{a_2}\delta_{a_3=0}$. Finally we take a suitable combination of this over different values of $l,m$ in order to get a solution of (11). To find the coefficients we need to write $\frac12(-1)^{a_1-1}\delta_{a_1\not=0}\delta_{a_2=0}$ as a combination of $\lambda^{a_1}\mu^{a_2}=\omega^{la_1+ma_2}$; the Fourier inversion formula gives
\[
\frac1{N^2}\sum\limits_{a_1=0}^{N-1}\sum\limits_{a_2=0}^{N-1}\frac12(-1)^{a_1-1}\delta_{a_1\not=0}\delta_{a_2=0}\omega^{-la_1-ma_2}
=\frac1{2N^2}\frac{1-\omega^l}{1+\omega^l}=\frac{-i}{2N^2}\tan\frac{\pi{}l}{N}
\]

So finally $*\pi((x)_0)$ is given by (8) where $r(g)$ is given by (5), $c$ by $(10)$ and $f=\d^*p$ with $p$ given by 
\[
p_a=\sum\limits_{l=1}^{N-1}\sum\limits_{m=0}^{N-1}\frac1{4N^2}\frac{1-\omega^l}{1+\omega^l}\frac{\omega^{la_1+ma_2}}{\nu^2-\nu^{-2}}
\left(\frac{\nu^{a_3}}{1-\nu^{N}}+\frac{(-\nu)^{a_3}}{1+\nu^N}
-\frac{\nu^{-a_3}}{1-\nu^{-N}}-\frac{(-\nu)^{-a_3}}{1+\nu^{-N}}
\right)
\]
where $\nu=\nu_{l,m}$ is given by (14). Rearranging,
\[
p_a=\frac1{4N^2}\sum\limits_{l=1}^{N-1}\sum\limits_{m=0}^{N-1}\tan\frac{\pi{}l}{N}\sin\frac{2\pi(la_1+ma_2)}{N}\frac1{\nu^2-\nu^{-2}}
\left(\frac{\nu^{a_3}+\nu^{N-a_3}}{1-\nu^{N}}+\frac{(-\nu)^{a_3}+(-\nu)^{N-a_3}}{1+\nu^N}
\right)
\]
\section{Invariant quantities}
According to the general theory of fluid algebras (see \S1), there are two invariant quantities, namely energy and helicity
\[(X,X)=\#(X\tm*X),\qquad(X,DX)=\#(X\tm*DX)=\#(X\tm\d{}X)\]
Explicitly,
\begin{align*}
(X,X)&=\#(X\tm*X)=\sum\limits_{a,b\in\Lambda}(u_a^{yz}u_b^{yz}+u_a^{zx}u_b^{zx}+u_a^{xy}u_b^{xy})2^{-\|a-b\|_1}\delta_{\|a-b\|_\infty\leq{}h}\\
&=\sum\limits_{a\in\Lambda}(u_a^{yz})^2+(u_a^{zx})^2+(u_a^{xy})^2+u_a^{yz}\bigg(u_{a+h\bfi}^{yz}+u_{a+h\bfj}^{yz}+u_{a+h\bfk}^{yz}+\frac12u_{a+h\bfi+h\bfj}^{yz}\\
&\quad+\frac12u_{a+h\bfi-h\bfj}^{yz}+\frac12u_{a+h\bfi+h\bfk}^{yz}+\frac12u_{a+h\bfi-h\bfk}^{yz}+\frac12u_{a+h\bfj+h\bfk}^{yz}+\frac12u_{a+h\bfj-h\bfk}^{yz}\\
&\quad+\frac14u_{a+h\bfi+h\bfj+h\bfk}^{yz}+\frac14u_{a+h\bfi+h\bfj-h\bfk}^{yz}+\frac14u_{a+h\bfi-h\bfj+h\bfk}^{yz}+\frac14u_{a+h\bfi-h\bfj-h\bfk}^{yz}\bigg)\\
&\quad+u_a^{zx}\bigg(u_{a+h\bfi}^{zx}+u_{a+h\bfj}^{zx}+u_{a+h\bfk}^{zx}+\frac12u_{a+h\bfi+h\bfj}^{zx}
+\frac12u_{a+h\bfi-h\bfj}^{zx}\\
&\quad+\frac12u_{a+h\bfi+h\bfk}^{zx}+\frac12u_{a+h\bfi-h\bfk}^{zx}+\frac12u_{a+h\bfj+h\bfk}^{zx}+\frac12u_{a+h\bfj-h\bfk}^{zx}\\
&\quad+\frac14u_{a+h\bfi+h\bfj+h\bfk}^{zx}+\frac14u_{a+h\bfi+h\bfj-h\bfk}^{zx}+\frac14u_{a+h\bfi-h\bfj+h\bfk}^{zx}+\frac14u_{a+h\bfi-h\bfj-h\bfk}^{zx}\bigg)\\
&\quad+u_a^{xy}\bigg(u_{a+h\bfi}^{xy}+u_{a+h\bfj}^{xy}+u_{a+h\bfk}^{xy}+\frac12u_{a+h\bfi+h\bfj}^{xy}
+\frac12u_{a+h\bfi-h\bfj}^{xy}\\
&\quad+\frac12u_{a+h\bfi+h\bfk}^{xy}+\frac12u_{a+h\bfi-h\bfk}^{xy}+\frac12u_{a+h\bfj+h\bfk}^{xy}+\frac12u_{a+h\bfj-h\bfk}^{xy}\\
&\quad+\frac14u_{a+h\bfi+h\bfj+h\bfk}^{xy}+\frac14u_{a+h\bfi+h\bfj-h\bfk}^{xy}+\frac14u_{a+h\bfi-h\bfj+h\bfk}^{xy}+\frac14u_{a+h\bfi-h\bfj-h\bfk}^{xy}\bigg)
\end{align*}

\section{Remarks on numerical simulations}

Coding of the fluid algebra was implemented by two independent groups using Runge--Kutta 4 algorithm for the time-stepping. In both cases, contrary to the theoretical calculations, the energy of the system blew up. Such is common in the case when an ODE is stiff, and often the problem is solved by implementing an implicit time-stepping. \\

An observation was made that if the points resulted by triple intersection is distinguished by the types of intersection, then the augmentation map can be changed without destroying the non-degeneracy of the inner product, which in turn would yield a more well-behaved ODE. \\

While this paper was being written, some of the authors developed another transverse intersection algebra TIA that is a differential graded algebra, though it is infinite dimensional and does not satisfy the crumbling property in \cite{ALS}. The fluid algebra associated with TIA seems to be more promising in the aspect of numerical stability. Coding of this new system will also be investigated as well. 

\section{Appendix by Dennis Sullivan}

{\bf Remark:A} There is a linear chain mapping forgetting the decoration from the new TIA to the  transverse intersection algebra EC (constructed in \cite{ALS}) which EC is a finite-dimensional commutative and associative algebra satisfying the product rule on all pairs of elements from the original complex, being TIA minus the decoration and the ideal elements. The EC algebra structure is the precursor of the algebra structure on TIA. This evolution was needed to improve the product rule and to  try to enable more stable fluid algebra computations. 

\

\noindent {\bf Remark:B} Those computations based on EC showed an instability in energy even though the system was mathematically conservative. There were two likely suspects for this instability in those computations: the odd subdivision  (introduced to make the inner product  of the fluid algebra (see \cite{S10}) nondegenerate and a dangerous structure constant in the  EC algebra venturing near a pole.
 The first can be eliminated by doing even subdivisions because in TIA the inner product is essentially non degenerate for even subdivisions.  The second suspect is  buffered away from the pole in TIA.  All of this in even period decompositions for   fluid algebra computations with the TIA  discretization; and these  will be made  when  the coding of TIA is completed.

\end{document}